\documentclass[12pt]{amsart}
\usepackage[T2A]{fontenc}
\usepackage[cp1251]{inputenc}
\usepackage[english]{babel}
\usepackage{amsmath,latexsym,amsthm,amsfonts,tipa,upgreek}
\usepackage{amssymb,amscd,graphpap, stmaryrd}

\theoremstyle{definition}

\begin{document}

\title{Recent progress in Subset Combinatorics of Groups}
\author{Igor Protasov, Ksenia Protasova}
%\subjclass{\bf MSC :54E15, 20F69}
\date{}
\address{Department of Computer Science and Cybernetics, Kyiv University, Volodymyrska 64, 01033, Kyiv, Ukraine}
\email{i.v.protasov@gmail.com; }
\address{Department of Computer Science and Cybernetics, Kyiv University, Volodymyrska 64, 01033, Kyiv, Ukraine}
\email{ksuha@freenet.com.ua}
\maketitle

\begin{abstract}
We systematize and analyze some results obtained in Subset Combinatorics of $G$ groups after publications the previous surveys [1-4]. The main topics: the dynamical and descriptive characterizations of subsets of a group relatively their combinatorial size, Ramsey-product subsets in connection with some general concept of recurrence in $G$-spaces, new ideals in the Boolean algebra $\mathcal{P}_{G}$ of all subsets of a group $G$  and in the Stone-$\check{C}$ech compactification $\beta G$  of $G$ , the combinatorial derivation.
\vspace{3 mm}

\vspace{3 mm}

\end{abstract}
\vspace{3 mm}

\section{Introduction}

In this paper, we systematize and analyze some results obtained in {\it Subset Combinatorics of Groups } after publications the surveys [1], [2], [3], [4]. The main topics: the descriptive and dynamical characterizations of subsets of a group with respect to their combinatorial size, Ramsey-product  subsets in connection with some general concept of recurrence, new ideals in the  Boolean algebra  $\mathcal{P}_{G}$ of all  subsets of $G$  and in the Stone-$\check{C}$ech compactification  $ \beta G$ of $G$, the combinatorial derivation.

In these investigations, the principal part play ultrafilters on a group $G$. On one hand, ultrafilters are using as a tool to get some purely combinatorial results. On the other hand, the  {\it Subset Combinatorics of Groups } allows to prove new facts about ultrafilters, in particular, about the Stone-$\check{C}$ech compactification  $\beta  G$ of $G$. In this connection, we recall some basic definitions concerning ultrafilters.

A {\it filter} $\mathcal{F}$ on a set $X$ is a family of subsets of $X$ such that
\vskip 5pt

$\bullet$  $\emptyset \notin \mathcal{F}$,  $X\in \mathcal{F}$;
\vskip 5pt

$\bullet$   $A, B \in \mathcal{F} \Longrightarrow A\bigcap B \in\mathcal{F}$;
\vskip 5pt

$\bullet$   $A \in \mathcal{F}$,  $A\subseteq C \Longrightarrow C\in\mathcal{F}$.
\vskip 7pt

The family of all filters on $X$  is partially ordered by inclusion.
A filter maximal in this ordering is called an {\it ultrafilter}. A filter $\mathcal{F}$ is an ultrafilter if and only if $X= A \bigcup B$  implies $A \in \mathcal{F}$ or $B\in \mathcal{F}$.

Now we endow $X$ with the discrete topology and identity the Stone-$\check{C}$ech compactification  $\beta X$ with the set of all ultrafilters on $X$. An ultrafilter  $\mathcal{F}$ is principal if there exists $x\in X$  such that $\mathcal{F}=\{A\subseteq  X:  x\in A\}$. Otherwise, $\bigcap\mathcal{F}=\emptyset$  and  $\mathcal{F}$ is called free. Thus, $X$   is identified with the set of all principal ultrafilters and the set of all free ultrafilter on $X$ is denoted by $X^{\ast}$.

To describe the topology on $\beta X$, given any $A \subseteq X$ we denote $\bar{A} = \{\mathcal{F}\in  X: A\in  \mathcal{F}\}$.  Then the set $\{ \bar{A}: A\subseteq X\}$
  is a base for the topology of $X$.  The characteristic topological property of $\beta X$: every mapping $f: X\longrightarrow  K$, $K$  is a compact Hausdorff space, can be extended to the  continuous mapping $f^{\beta}: \beta X\longrightarrow  K$.

Given a filter $\varphi$ on $X$, the set $\bar{\varphi} = \{p\in\beta X: \varphi\subseteq p\}$
 is closed in $\beta X$,  and for every non-empty closed subset $K$ of $\beta X$,  there is a filter $\varphi$ on $X$  such that $\bar{\varphi}=K$.

Now let $G$  be a discrete group. Using the characteristic property of $\beta G$, we can extend the group multiplication on $G$  to the semigroup multiplication on $\beta G$  in such a way that, for every $g\in G$, the mapping $\beta G\longrightarrow G:  p\longmapsto gp$  is  continuous and, for every $q\in  \beta G$, the
mapping  $\beta G \longrightarrow \beta G:  p\longmapsto pq$
is continuous.

To define  the product $pq$ of ultrafilters $p$  and $q$, we take an arbitrary $P\in p$  and, for each $x\in P$, pick some $Q _{x}\in q$. Then, $\bigcup_{x\in P}  x Q_{x}$ is a member of $pq$, and each member of $pq$  contains some subsets of this form.

For properties of the compact right topological semigroup $\beta G$ and a plenty of its combinatorial application see [5].

%zzzzzzzzzzzzzzzzzzzzzzzzzzzzzzzzzzzzz Секция 2

\section{Diversity of subsets and ultracompanions}

Let $G$  be a group with the identity $e$, $\mathcal{F}_{G}$ denotes the family of all finite subsets of $G$.
 We say that a subset $A$ of $G$ is \vskip 5pt

 \begin{itemize}
\item{} {\it large} if  $G=FA$  for some $F\in \mathcal{F}_G$;\vskip 5pt

\item{} {\it small} if  $L\setminus A$ is  large for every large subset $L$ ;
\vskip 5pt

\item{} {\it extralarge} if  $G\setminus  A$ is small;\vskip 5pt

\item{} {\it thin} if  $gA\cap A$ is  finite for  each $g\in G\setminus \{e\}$;\vskip 5pt

\item{} {\it thick} if,  for every  $F\in \mathcal{F}_G$, there exists $a\in A$ such that  $Fa \subseteq A$;\vskip 5pt

\item{} {\it prethick} if  $FA$ is thick for some $F\in \mathcal{F}_G$;
\vskip 5pt

\item{} {\it $n$-thin}, $n\in {\mathbb N}$
 if, for every   distinct elements   $g_0 , \dots , g_n \in  G$,
the set  $g_0  A \cap \dots \cap g_n  A$    is  finite; \vskip 5pt

\item{} {\it sparse} if,  for every   infinite subset   $X$ of $G$, there exists a  finite subset $F\subset X$  such that $\bigcap _{g\in F}  gA$ is  finite. \vskip 5pt
\end{itemize}
\vskip 10pt

{\bf Remark 2.1.} In {\it Topological dynamics}, large subsets are known as syndetic, and a subset is small if and only if it fails to be piecewise syndetic. In \cite{b4}, the authors use the dynamical terminology.
\vskip 6pt

All above definitions can be unified with usage the following notion \cite{b6}. Given a subset $A$ of a group $G$ and an ultrafilter $p\in G^*$, we define a $p$-{\it companion} of $A$  by
$$ \Delta _p (A)= A^* \cap Gp=\{gp: g\in G,  \   A\in gp  \}.  $$
\vskip 5pt

Then, for every infinite group $G$, the following statement hold:
\vskip 10pt

 \begin{itemize}
\item{} $A$ is large if  and only if  $\Delta _p (A) \neq  \emptyset  $ for each    $p\in G^*$;\vskip 5pt

\item{} $A$ is  small if   and only if,
for every  $p\in   G^*$  and every    $F\in  \mathcal{F}_G$,  we have $\Delta _p (F A) \neq  Gp  $;\vskip 5pt

\item{} $A$ is thick if and only if, there exist $p\in G^{\ast}$ such that $\Delta _p (A) = Gp $;  \vskip 5pt

\item{} $A$ is thin if and only if, $\Delta _p (A) \leq 1 $  for every $p\in G^*$;\vskip 5pt

\item{} $A$ is $n$-thin if and only if,  $\Delta _p (A) \leq n $   for every $p\in G^*$; \vskip 5pt

\item{} $A$ is sparse if and only if, $\Delta _p (A) $ is finite for each $p\in G^*$.\vskip 5pt

\end{itemize}

\vskip 10pt

 Following \cite{b1}, we say that a subset $A$  of $G$ is {\it scattered} if, for every infinite subset $X$ of $A$, there is $p\in X^* $  such that $\Delta _p (X)$ is finite. Equivalently [7, Theorem 1],  $A$ is scattered if each subset  $\Delta _p (A)$ is discrete in $G ^* $.

\vskip 10pt

{\it Comments}.  For motivations of above definitions see  \cite{b1}, for more delicate classification of subsets of a group and $G$-spaces see \cite{b2}, \cite{b8}.

%zzzzzzzzzzzzzzzzzzzzzzzzzzzzzzzzzzzzz Секция 3

\section{The descriptive look at the size of subsets of groups}

Given a group $G$, we denote by ${\bf P}_{G}$ and ${\bf F}_{G}$ the Boolean algebra of all subsets of $G$ and its ideal of all finite subsets.  We endow ${\bf P}_{G}$ with the topology arising from identification (via characteristic functions) of ${\bf P}_{G}$ with $\{0,1\}^{G}$. For $K\in F_{G}$ the sets   $$ \{X\in {\bf P}_{G}: K\subseteq X\},  \  \  \{X\in {\bf P}_{G}: X\cap K=\emptyset\}$$  form the subbase of this topology.

After the topologization, each family $\mathcal{F}$ of subsets of a group  $G$ can be considered as a subspace of  ${\bf P}_{G}$,  so one can ask about the Borel complexity of $\mathcal{F}$, the question typical in the {\it Descriptive Set Theory} (see \cite{b9}). We ask these questions for the most intensively studied  families in  {\it Combinatorics of Groups.}

For a group $G$, we denote by ${\bf L}_{G}$, ${\bf EL}_{G}$, ${\bf S}_{G}$,  ${\bf T}_{G}$, ${\bf PT}_{G}$  the sets of all large, extralarge, small, thick and prethick subsets of $G$, respectively.

\vskip 10pt

{\bf Theorem 3.1.}{\it For a countable group $G$, we have: ${\bf L}_{G}$ is  $F_{\sigma}$, ${\bf T}_{G}$ is $G_{\delta}$, ${\bf PT}_{G}$  is $G_{\delta\sigma}$, ${\bf S}_{G}$ and ${\bf EL}_{G}$  are} $F_{\sigma\delta}$.
%$\Box$
\vskip 10pt

A subset $A$ of a group $G$ is called
\begin{itemize}
\item {\em  $P$-small\/}  if there exists  an injective sequence   $(g_{n})_{n\in\omega}$  in $G$  such  that the subsets $\{g_{n} A : n\in\omega \}$ are pairwise disjoint;
\item {\em weakly $P$-small\/}  if, for any $n\in\omega$,  there exists $g_{0},\ldots , g_{n}$ such that the subsets $g_{0}A,\ldots , g_{n}A$ are pairwise disjoint;
\item {\em almost $P$-small\/}  if there exists  an injective sequence  $(g_{n})_{n\in\omega} $ in $G$  such that  $g_{n}A\cap g_{m}A$   is finite for all distinct $n,m$;
\item {\em near $P$-small\/}  if,  for every $n\in\omega$, there exists $g_{0},\ldots , g_{n}$ such that $g_{i}A\cap g_{j}A$ is finite for all distinct  $i,j\in \{0,\ldots,n\}$.
\end{itemize}

\vskip 10pt

Every infinite group $G$  contains a weakly  $P$-small set, which is not $P$-small, see   \cite{b10}. Each almost $P$-small subset  can be partitioned into two  $P$-small  subsets \cite{b8}. Every countable Abelian group contains a near  $P$-small subset which is neither weakly nor almost $P$-small \cite{b11}.
\vskip 10pt

{\bf Theorem 3.2.}{\it For a countable group $G$,  the sets of thin, weakly $P$-small  and near $P$-small subsets of $G$ are $F_{\delta\sigma}$.}

\vspace{5 mm}

We recall that a topological space $X$ is {\it Polish} if $X$ is homeomorphic to a separable complete metric space.
A subset $A$ of a   topological space $X$ is {\it analytic} if $A$ is a continuous image of some Polish space, and $A$ is {\it coanalytic} if $X\setminus A$ is analytic.

Using the classical tree technique \cite{b9} adopted to groups in \cite{b12}, we get.\vspace{3 mm}
\vskip 10pt

{\bf Theorem 3.3.} {\it For a countable group $G$, the ideal of sparse subsets  is coanalytic and the set of $P$-small  subsets is analytic in   ${\bf P}_{G}$.}
\vspace{5 mm}

Given a discrete group $G$, we identify the Stone-$\check{C}$ech compactification  $\beta G$ with the set of all ultrafilters on $G$ and consider  $\beta G$ as a right-topological semigroup (see  Introduction). Each non-empty  closed subspace $X$ of $\beta G$ is determined by some filter $\varphi$ on $G$: $$X=\bigcap\{\overline{\Phi} : \Phi\in\varphi\},  \  \   \overline{\Phi}=\{p\in\beta G: \Phi\in p \}. $$

On the other hand, each filter $\varphi$ on $G$ is a subspace of
${\bf P}_{G}$, so we can ask about complexity of $X$ as the complexity   of
 $\varphi$ in ${\bf P}_{G}$.

The semigroup $\beta G$  has the minimal ideal $K_{G}$ which play one of the key parts in combinatorial  applications of   $\beta G$. By \cite{b5}, Theorem 1.5, the closure $cl (K_{G})$ is determined by the filter of all extralarge  subsets of $G$.  If $G$ is countable, applying Theorem 3.1, we conclude that  $cl (K_{G})$  has the Borel complexity  $F_{\sigma\delta}$.

An ultrafilter $p$ on $G$ is called {\it strongly  prime} if $p\notin  cl(G^{\ast} G^{\ast})$, where   $G^{\ast}$ is a semigroup of all free ultrafilters on $G$.  We put  $X=  cl(G^{\ast} G^{\ast})$  and choose the filter  $\varphi_{X}$ which determine $X$. By \cite{b13}, $A\in \varphi_{X}$ if and only if $G\backslash A$ is sparse. If $G$  is countable,  applying Theorem 3.3, we conclude  that $\varphi_{X}$  is coanalitic in  ${\bf P}_{G}$.

Let  $(g_{n})_{n\in\omega}$ be an injective sequence in $G$. The set $$\{g_{i_{1}} g_{i_{2}}\ldots  g_{i_{n}}: 0 \leq i_{1}< i_{2}< \ldots < i_{n}<\omega \}$$
is called an  {\it FP-set}. By the Hindman Theorem 5.8 \cite{b5}, for every finite partition of $G$, at least one cell of the partition contains an $FP$-set. We denote
by   ${\bf FP}_{G}$ the family of all subsets of $G$ containing some $FP$-set. A subset $A$ of $G$  belongs to  ${\bf FP}_{G}$ if  and only if $A$ is an element of some idempotent of  $\beta G$. By analogy with Theorem 3.3, we can prove that  ${\bf FP}_{G}$ is analytic in ${\bf P}_{G}$.
\vskip 10pt

{\it Comments.} This section reflects the results from \cite{b14}.

%zzzzzzzzzzzzzzzzzzzzzzzzzzzzzzzzzzzzz Секция 4

\section{The dynamical look at the subsets of a group}

Let $G$ be a group. A topological space $X$ is called a $G$-{\it space} if there is the
action $X \times G \longrightarrow X: (x, g) \longmapsto xg$ such that, for each $g\in G$, the mapping
$X \longrightarrow X: x \longmapsto xg$ is continuous.

Given any $x\in X$ and $U \subseteq X$, we set
$$[U]_{x} = \{g \in G: xg\in U\}$$
and denote
$$O(x) = \{xg : g \in G\}, T (x) = cl O(x),$$

$W(x) = \{y \in T (X) : [U]_{x}$    is infinite for each neighbourhood $U$ of  $y\}$.

We recall also that $x \in X$ is a {\it recurrent point} if $x \in W(x)$.

Now we  identify $\mathcal{P}_G$ with the space $\{0, 1\}^{G}$, endow
$\mathcal{P}_G$ with the product topology and  consider $\mathcal{P}_G$ as a $G$-space with the action defined by
$$A \longmapsto Ag, \  Ag = \{ ag : a \in A\}.$$
We say that a subset $A$ of $G$ is {\it recurrent} if  $A$ is a recurrent point in $(\mathcal{P}_G, G)$.

All groups in this sections are supposed to be infinite.
\vskip 10pt

{\bf Theorem 4.1.} {\it For a subset $A$ of a group $G$, the following statements hold
\vspace{5 mm}

(i)  $A$ is finite if and only if  $W(A) =\emptyset $;
\vspace{5 mm}

(ii) $A$ is thick if and only if  $G \in W(A)$.}
\vskip 10pt

{\bf Theorem 4.2.} {\it  For a subset $A$ of a group $G$, the following statements hold
\vspace{5 mm}

(i) $A$ is $n$-thin if and only if $|Y | \leq n$ for every $Y \in W(A)$;
\vspace{5 mm}

(ii) $A$ is sparse if and only if each subset $Y \in W(A)$ is finite;
\vspace{5 mm}

(iii) $A$ is scattered if and only if, for every subset $B \subseteq A$ there exists $Y \in \mathcal{F} _{G}$ in the closure of}  $\{Bb^{-1} : b \in B \}.$
\vspace{5 mm}

Let $(g_{n})_{n\in\omega}$ be an injective sequence in $G$.
The set
$$FP(g_{n})_{n\in\omega} = \{g_{i_{1}}g_{i_{2}} \ldots g_{i_{n}} : 0 \leq i_{1} < i_{2} < \ldots < i_{n} < \omega\}$$
is called an $FP$-set.

Given a sequence $(b_{n})_{n\in\omega}$  in $G$, the set
$$\{g_{i_{1}}g_{i_{2}} \ldots g_{i_{n}} b_{i_{n}} : 0 \leq i_{1} < i_{2} < \ldots < i_{n} < \omega\}$$
is called a {\it  (right) piecewise shifted} $FP$-set [7].
\vskip 10pt

{\bf Theorem 4.3.}
{\it For a subset $A$ of a group $G$, the following statements hold
\vspace{5 mm}

(i) $A$ is not $n$-thin if and only if there exist $F \in [G]^{n+1}$  and an injective
sequence $(x_{n})_{n<\omega}$ in $G$ such that $Fx_{n} \subseteq A$ for each $n \in  \omega$;
\vspace{5 mm}

(ii) $A$ is not sparse if and only if there exists two injective sequences $(x_{n})_{n<\omega}$
and $(y_{n})_{n<\omega}$ such that $x_{n}y_{m} \in A$ for each $0 \leq n \leq m < \omega$;
\vspace{5 mm}

(iii) $A$ is not scattered if and only if $A$ contains a piecewise shifted $FP$-set;
\vspace{5 mm}

(iv) $A$ contains a recurrent subset if and only if there exists $x \in A$ and an
$FP$-set $Y$ such that $xY \subseteq A$.}
\vskip 10pt

{\bf Corollary 4.1. }
{\it  Every scattered subset of a group G has no recurrent points.}
\vspace{5 mm}

{\bf  Remark 4.1.} By [4, Theorem 2], every scattered subset $A$ of an amenable group
$G$ is absolute null, i.e. $\mu(A) = 0$ for every left invariant Banach measure $\mu$ on
$G$. But this statement could not be generalized to subsets with no recurrent
points. By [17, Theorem 11.6], there is a subset $A$ of $\mathbb{Z}$ of positive Banach
measure such that $(a + B) \setminus A \neq \emptyset $ for any $FP$-set $B$. By Theorem 4.3(iv), $A$
has no recurrent subsets.

\vskip 10pt

{\bf  Remark 4.2.} Let $G$ be an arbitrary infinite group. In \cite{b15}, we
 constructed two injective
sequences $(x_{n})_{n\in\omega}$, $(y_{n})_{n\in\omega}$ in $G$ such the set $\{x_{n}y_{m} : 0 \leq n \leq m < \omega\}$ is scattered. By Theorem 4.3(ii), this subset is not sparse.
\vspace{5 mm}

{\it Comments}. This section reflects the first part of \cite{b15}.

%zzzzzzzzzzzzzzzzzzzzzzzzzzzzzzzzzzzz Секция 5

\section{Ramsey-product subsets and recurrence}

In this section, all groups under consideration are supposed to be infinite; a countable set means a    countably infinite set.

Let $G$ be a group and let $\overrightarrow{m}= (m_{1} \ldots, m_{k}) \in\mathbb{Z}^{k}$ be a number vector of length
$k\in \mathbb{N}$.
We say that a subset $A$ of a group $G$ is a {\it Ramsey $\overrightarrow{m}$-product subset}
if  every infinite subset $X$ of  $G$ contains pairwise
distinct elements    $x_{1},\ldots, x_{k} \in X$ such that,
  $$x^{m_{1}} _{\sigma(1)} \  x^{m_{2}} _{\sigma(2)}  \ldots  x^{m_{k}} _{\sigma(k)} \in A$$
for every substitution $\sigma\in S_{k}$.

 \vskip 10pt

{\bf  Theorem 5.1.}. {\it  For a group $G$  and a number vector $\overrightarrow{m}=(m_{1},\ldots , m_{k} )  \in\mathbb{Z}^{k}$, the following statements hold:
\vskip 15pt
$(i)$  a subset $A$  of  $G$ is a  Ramsey $\overrightarrow{m}$-product subset if  and only  if every infinite subset $X$ of $G$ contains a countable subset $Y$ such that $y_{1}^{m_{1}} \ldots  y_{k}^{m_{k}}\in A$ for any distinct elements
$y_{1}, \ldots  , y_{k} \in Y$.
\vskip 15pt

$(ii)$  the family  $ \varphi _{\overrightarrow{m}}$  of  all Ramsey $\overrightarrow{m}$-product subsets of $G$ is a filter. }\vskip 15pt

For $t\in\mathbb{Z}$ and $q\in G^{\ast}$ we denote by $q^{\wedge}t$ the ultrafilter with the base $\{x^{t}: x\in Q\}$,  $ \ Q\in q$.
  Warning: $q^{\wedge}t$ and $q^{t}$ are different things. Certainly,  $q^{\wedge}t=q^{t}$ only if $t\in\{-1,0,1\}$.

We remind the reader that, for a filter $\varphi$ on $G$,  $\overline{\varphi}= \{p\in\beta G:\varphi\subseteq p\}$.

  {\bf  Theorem 5.2.} {\it  For every group $G$ and any number vector $\overrightarrow{m}=(m_{1}, \ldots , m_{k}) \in\mathbb{Z}^{k}$,  we have}
$$\overline{\varphi}_{\overrightarrow{m}} \ = \ cl\{(q^{\wedge}m_{1}) \  \ldots \  (q^{\wedge}m_{k}):  \  q\in \  G^{\ast}\}. $$
 \vskip 5pt

Now we consider some special cases of vectors $\vec{m}$.
\vskip 10pt

{\bf Proposition 5.1.} {\it For any totally bounded topological group G, any neighborhood $U$  of
the identity $e$ of $G$ is a Ramsey
 $\vec{m}$-product subset for any  vector $\vec{m} = (m_{1}, \ldots , m_{k})$
 such that}
 $ m_{1} + \ldots + m_{k} = 0.$
\vskip 5pt

We recall that a {\it quasi-topological group} is a group $G$
endowed with a topology such that, for any $a, b \in G$ and $\varepsilon \in \{-1, 1\}$,
 the mapping $G \longrightarrow G:
x \longmapsto ax^{\varepsilon}b$, is continuous.
\vskip 10pt

{\bf Proposition 5.2.} {\it  The closure $\bar{A}$ of any Ramsey $(-1, 1)$-product set $A$ in a quasi-topological
group $G$ is a neighborhood of the identity.}
\vskip 10pt

{\bf Proposition 5.3.} {\it  Let $\vec{m} = (m_{1}, \ldots , m_{k})$
 be a number vector and $s = m_{1} + \ldots + m_{k}$.
For any Ramsey  $\vec{m}$-product subset A of a group G, the set $\{x^{s} : x \in G\}$ is contained in the
closure of  $A$ in
any non-discrete group topology on $G$.}

\vskip 10pt

{\bf Proposition 5.4.} {\it  Let $G$  be the Boolean group of all finite subsets of  $\mathbb{Z}$, endowed with the group operation of symmetric difference.  The set
$$A = G\setminus \{ \{x, y\} : x, y \in \mathbb{Z}, 0 \neq x - y \in \{z^{3} : z \in \mathbb{Z}\}\}$$
has the following properties:
\vskip 5pt

(i)  $A$ is a Ramsey $\vec{m}$-product  for any vector $\vec{m}=(m_{1}, \ldots ,m_{k}) \in (2\mathbb{Z} + 1)^{k}$ of length $k \geq 2$;
\vskip 5pt

(ii)  $A$  does not contain the difference $BB ^{-1}$ of any large subset $B$ of  $G$;
\vskip 5pt

(iii)  $A$  is not a neighborhood of zero in a totally bounded group topology on} $G$.
\vskip 10pt

Now we show how Ramsey  $(-1, 1)$-product sets arise in some  general concept of recurrence on $G$-spaces.

Let $G$ be a group with the identity $e$  and let $X$ be a $G$-space with the action
$G\times X\longrightarrow X$, $(g,x)\longmapsto gx$. If $X=G$  and $gx$ is the product of $g$ and $x$ then $X$ is called a {\it left regular $G$-space}.

Given a $G$-space $X$, a family $\mathfrak{F}$ of subset of $X$ and $A\in \mathfrak{F}$, we denote
\begin{eqnarray}
\nonumber \Delta_{\mathfrak{F}}(A)=\{g\in G: gB\subseteq A  \text{  for  \ some }  B\in\mathfrak{F}, B\subseteq A\}.
\end{eqnarray}

Clearly,  $e\in \Delta_{\mathfrak{F}}(A)$ and if $\mathfrak{F}$ is upward directed $(A\in \mathfrak{F}$, $A\subseteq C$ imply $C\in \mathfrak{F})$  and if $\mathfrak{F}$ is $G$-invariant  $(A\in \mathfrak{F}$, $g\in G$ imply $gA\in \mathfrak{F})$ then
\begin{eqnarray}
\nonumber \Delta_{\mathfrak{F}}(A)=\{g\in G: gA\cap A \in \mathfrak{F}\},  \Delta_{\mathfrak{F}}(A)= (\Delta_{\mathfrak{F}}(A))^{-1}.
\end{eqnarray}

If $X$ is a left regular $G$-space and $\emptyset\notin  \mathfrak{F}$ then $\Delta_{\mathfrak{F}}(A)\subseteq A A^{-1}.$

For a $G$-space $X$ and a family $\mathfrak{F}$ of subsets of $X$, we say that a subset $R$ of $G$ is $\mathfrak{F}$-{\it recurrent} if $\Delta_{\mathfrak{F}}(A)\cap R \neq \emptyset$  for every $A\in \mathfrak{F}$. We denote by $\mathfrak{R}_{\mathfrak{F}}$  the filter on $G$ with the base $\cap\{\Delta_{\mathfrak{F}\prime}(A):  A\in \mathfrak{F}^{\prime}\},$  where $\mathfrak{F}^{\prime}$ is a finite subfamily of $\mathfrak{F}$, and note that, for an ultrafilter $p$ on $G$,  $\mathfrak{R}_{\mathfrak{F}}\in p$ if and only if each member of $p$ is $\mathfrak{F}$-recurrent.

The notion of an $\mathfrak{F}$-recurrent subset is well-known in the case in which $G$ is an amenable group, $X$ is a left regular $G$-space  and $\mathfrak{F}=\{ A\subseteq X: \mu(A)>0$  for some left invariant Banach measure $\mu$ on $X\}$. See \cite{b16},  \cite{b17}, \cite{b18} for historical background.

We recall \cite{b19} that a filter $\varphi$ on a group $G$  is {\it left topological}  if $\varphi$ is a base at the identity $e$ for some (uniquely  defined) left translation invariant (each left shift $x\longmapsto gx$ is continuous) topology on $G$. If $\varphi$ is left topological then $\overline{\varphi}$ is a subsemigroup of $\beta G$  \cite{b19}.
If $G=X$  and a filter $\varphi$ is left topological then $\varphi= \mathfrak{R}_{\varphi}$.
\vskip 10pt

{\bf Proposition 5.5.} {\it For every $G$-space $X$ and any family $\mathfrak{F}$ of subsets of $X$, the filter $\mathfrak{R}_{\mathfrak{F}}$ is left topological. }
\vskip 5pt

Let $X$ be a $G$-space  and let $\mathfrak{F}$ be a family of subsets of $X$. We say that a family $\mathfrak{F}^{\prime}$ of subsets of $X$ is $\mathfrak{F}$-{\it disjoint} if $A\cap B\notin \mathfrak{F}$  for any distinct $A,B \in \mathfrak{F}^{\prime}$.

A family $\mathfrak{F}^{\prime}$ of subsets of $X$ is called $\mathfrak{F}$-{\it packing large} if,  for each   $A\in\mathfrak{F}^{\prime}$, any $\mathfrak{F}$-disjoint family of subsets of $X$  of the form $gA,$ $g\in G$ is finite.
\vskip 10pt

{\bf Proposition 5.6.} {\it Let $X$  be a $G$-space and let $\mathfrak{F}$ be a $G$-invariant upward directed family of subsets of $X$. Then $\mathfrak{F}$ is $\mathfrak{F}$-packing large if and only if, for each $A\in \mathfrak{F}$, the set $ \ \triangle_{\mathfrak{F}}(A) \ $  is a  Ramsey (-1,1)-product set}.

\vskip 5pt

Applying Theorem 5.2, we conclude that  $ \ \triangle_{\mathfrak{F}}(A) \ $   contains all ultrafilters of the form $q^{-1}q$,  $q\in G^{\ast}$, and in the case $X=G$,  $G$  is amenable and $\mathfrak{F}$
 is the family of all subsets of positive Banach measure,  we get Theorem 3.14 from \cite{b18}.
\vskip 5pt

{\it Comments}. The proofs of all above statements can be find in \cite{b20},  \cite{b21}.

%zzzzzzzzzzzzzzzzzzzzzzzzzzzzzzzzzzzz Секция 6

\section{  Ideals in $\mathcal{P} _G$  and $\beta G$ }

We recall that a family $ \mathcal{I}$ of subsets of a set $X$ is an {\it ideal}  in the Boolean algebra $\mathcal{P} _G$ of all subsets of $G$ if $G \notin \mathcal{I}$
and $A\in \mathcal{I}$, $B\in \mathcal{I}$, $C\subseteq A$  imply $A\cup B\in \mathcal{I}$, $C\in \mathcal{I}$.  A family  $\varphi$ of subsets of $G$ is a  filter if and only if
 the family $\{ X\setminus  A: A \in \varphi \}$ is an ideal.

For an infinite group $G$, an ideal $\mathcal{I}$ in $ \mathcal{P}_G$ is called {\it left (right) translation invariant} if $gA\in \mathcal{I}$ ($Ag\in \mathcal{I}$) for all $g\in G$, $A\in \mathcal{I}$. If $\mathcal{I}$ is left and right translation invariant then $\mathcal{I}$ is called {\it translation invariant}. Clearly, each left (right) translation invariant ideal of $G$ contains the ideal  $\mathcal{F} _G$ of all finite subsets of $G$. An ideal $\mathcal{I}$ in $\mathcal{P} _G$ is called a {\it group ideal} if  $\mathcal{F} _G\subseteq \mathcal{I}$
and if $A\in \mathcal{I}$, $B\in \mathcal{I}$  then $AB ^{-1}\in \mathcal{I}$.

Now we endow $G$ with the discrete topology and
 use the standard extension  of the  multiplication on $G$ to the semigroup multiplication on  $\beta G$, see Introduction.

It follows directly from the definition of the multiplication in $\beta G$ that   $G^{*}$,   $ \overline{ G^{*} G^{*}}$ are ideals  in the semigroup $\beta G$, and  $G^{*}$  is the unique  maximal closed ideal in  $\beta G$. By Theorem 4.44 from [5],  the closure  $ \overline{ K(\beta G)}$  of the minimal ideal    $K(G)$  of   $\beta G$  is an ideal, so $ \overline{ K(\beta G)}$ is the smallest closed ideal in $\beta G$.  For the structure of    $ \overline{ K(\beta G)}$  and some other ideals in  $\beta G$  see [5, Sections 4,6].
\vskip 10pt

For an ideal $\mathcal{I}$ in $\mathcal{P}_G$, we put
$$\mathcal{I}^{\wedge} = \{p\in  \beta G: G\setminus A \in p  \text{  for each  }  A\in \mathcal{I} \}, $$
 and use the following observations: \vskip 5pt

\begin{itemize}
\item{} $\mathcal{I}$ is left translation invariant if and only if $\mathcal{I}^{\wedge}$    is a left ideal of the semigroup $\beta G$ ; \vskip 5pt

\item{}  $\mathcal{I}$ is right  translation invariant if and only if
$(\mathcal{I}^{\wedge})G\subseteq  \mathcal{I}^{\wedge}$.

\end{itemize}
\vskip 5pt

We use also the inverse to $^\wedge$  mapping  $^\vee$.  For a closed  subset $K$ of $\beta G$, we take
 the unique filter $\varphi$  on $G$  such that $K=\overline{\varphi}$  and put
$$K^{\vee} = \{G\setminus A : A\in    \varphi\}.$$

In this section, all groups under consideration are suppose to be infinite.

We denote by $Sm_G$, $Sc_G$, $Sp_G$ the families of all small, scattered and sparse subsets of a group $G$. These families are translation invariant ideals in $\mathcal{P}_G$
 (see [6, Proposition 1 ]),  and for every group $G$, the following inclusions are strict [6,  Proposition 12]
$$ Sp_G \subset  Sc_G  \subset  Sm_G . $$

We say that a subset $A$  of $G$ is
{\it finitely thin} if $A$  is $n$-thin for some $n\in \  \mathbb{N}$.
The family $ FT_G$ of  all  finitely  thin subsets of $G$  is a translation invariant ideal in
$\mathcal{P}_G$   which contains the ideal $<T_G>$ generated by the family of all thin subsets of $G$.
 By  [22, Theorem 1.2] and [23, Theorem 3], if $G$ is either countable or Abelian and $|G| < \aleph _\omega$ then
$FT_G = <T_G>$. By [23, Example 3], there exists an Abelian group $G$ of cardinality  $\aleph _\omega$
such that $<T_G>\subset FT_G $.

\vskip 10pt

% ZZZZZZZZZZZZZZ heorem 6.1.
{\bf Theorem 6.1.}
{\it For every group $G$, we have $Sm_G ^\wedge
= \overline{ K(\beta G)}$. }\vskip 10pt

This is Theorem 4.40 from \cite{b5}  in the form given in [24, Theorem 12.5].

\vskip 10pt
% ZZZZZZZZZZZZZZ heorem 2.2.
{\bf Theorem 6.2.}
{\it  For every group $G$,
 $Sp_G ^\wedge= \overline{ G ^*G ^*}$;}

\vskip 10pt

This is Theorem 10  from \cite{b13}.

%zzzzzzzzzzzzzzzzzzzzzzzzzzzzzzzzzzzzz

{\bf 6.1. Between $\overline{ G ^*G ^*}$  and  $G ^*$.}
\vskip 10pt

% ZZZZZZZZZZZZZZ heorem 3.1.
{\bf Theorem 6.3.}
{\it  For every group $G$, the following statements hold:
\vskip 5pt

$(i)$   if  $\mathcal{I}$ is a  left  translation invariant ideal in $\mathcal{P}_G$  and $\mathcal{I} \neq \mathcal{F}_G$ then there exists a left translation invariant ideal  $\mathcal{J}$ in  $\mathcal{P}_G$ such  that  $\mathcal{F}_G \subset \mathcal{J} \subset  \mathcal{I} $ and   $\mathcal{J} \subset  Sp_G $;
\vskip 5pt

$(ii)$  if  $\mathcal{I}$ is a  right  translation invariant ideal in $\mathcal{P}_G$  and $\mathcal{I} \neq \mathcal{F}_G$ then there exists a right  translation invariant  $\mathcal{J}$ in  $\mathcal{P}_G$ such  that  $\mathcal{F}_G \subset \mathcal{J} \subset  \mathcal{I} $;
\vskip 5pt

$(iii)$  if $G$  is  either countable or Abelian and
$\mathcal{I} $ is a  translation invariant ideal in $\mathcal{P}_G$ such  that
 $\mathcal{I} \neq \mathcal{F}_G$ then there exists a   translation invariant   ideal  $\mathcal{J}$ in  $\mathcal{P}_G$ such  that  $\mathcal{F}_G \subset \mathcal{J} \subset  \mathcal{I} $ and   $\mathcal{J} \subset  Sp_G $;}
\vskip 10pt

%zzzzzzzzzzzzzzzzzzzzzzTheorem 6.4.
{\bf Theorem 6.4.}
{\it  For every group $G$, the following statements hold:
\vskip 5pt

$(i)$   if  $L$ is a  closed left ideal in   $\beta G$ such that $L\subset G^*$ then there exists a closed left ideal
$L ^\prime$  of    $\beta G$   such that $L\subset L^\prime  \subset G^* $,   $\overline{ G ^*G ^*}\subset  L ^\prime$;

\vskip 5pt

$(ii)$   if  $R$ is a  closed subset of  $G^*$ such that $R\neq G^*$ and  $RG\subseteq  R$
 then there exists a closed subset
$R ^\prime$  of    $G^*$   such that
 $R\subset R^\prime  \subset G^* $,   $ R ^\prime G\subseteq  R $;

\vskip 5pt

$(iii)$   if  $G$ is  either countable or Abelian and  $I$  is
a  closed
ideal in   $\beta G$ such that
$I\subset G^*$ then there exists a closed ideal
$I ^\prime$  in    $\beta G$   such that $I\subset I^\prime  \subset G^* $,
$\overline{ G ^*G ^*}\subset  I $.}
\vskip 15pt

For a cardinal $\kappa$, $S_\kappa$ denotes the group of all permutations  of $\kappa$.
\vskip 15pt

%zzzzzzzzzzzzzzzzzzzzzzzzzzzzzzzzTheorem 3.3.
{\bf Theorem 6.5.}
{\it  For every infinite  cardinal $\kappa$,  there exists a closed ideal $I$ in $\beta S_\kappa$ such that
\vskip 10pt

$(i)$ $S_\kappa ^*  S_{\kappa}^*\subset I  $;
\vskip 10pt

$(ii)$ if $M$ is a closed ideal in $\beta S_\kappa $ and $I\subseteq M\subseteq  G^*$  then either $M=I$ or $M= S_\kappa ^*$.}
\vskip 10pt

%zzzzzzzzzzzzzzzzzzzzzTheorem
{\bf Theorem 6.6.}
{\it  For every group $G$, we have $ {\it FT_G \subset Sp _G}  $ so
$\overline{ G ^*G ^*}\subset   FT_G ^{\wedge}$.}
\vskip 10pt

For subsets $X, Y$  of a group $G$, we say that the product $XY$ is an $n$-{\it stripe} if $|X|=n$,
$n\in \mathbb{N}$  and $|Y|=\omega$.   It is easy to see that a subset   $A$  of $G$  is $n$-thin  if and only if $A$  has no $(n+1)$-stripes. Thus, $p\in FT_G ^\wedge$ is and only if each member $P\in p$  has an   $n$-stripe for every $n\in \mathbb{N}$.

We say that $XY$ is an $(n,m)$-{\it rectangle}  if $|X|=n$, $|Y|=m \ $, $ \ n,m\in \mathbb{N}$. We say that a subset $A$ of $G$ {\it has bounded rectangles} if there is $n\in \mathbb{N}$  such that $A$  has no $(n, n)$-rectangles (and so $(n, m)$-rectangles for each $m>n$).

We denote by $BR_G$ the family of all subsets of $G$ with bounded rectangles.

\vskip 15pt

%ZZZZZZZZZZZZZZZZZZZZZZZZ
{\bf Theorem 6.7.}
{\it  For a group $G$, the following statements hold:\vskip 10pt

$(i)$ $BR_G$ is a  translation invariant ideal  in $ \mathcal{P}_G$  ;

\vskip 10pt

$(ii)$  $BR_G^\wedge$ is a closed ideal in $\beta G$  and  $p\in BR_G^\wedge$ if and only if each member $P\in p$ has an $(n,n)$-rectangle for  every $n\in   \mathbb{N} $;

\vskip 10pt

$(iii)$   $BR_G \subset FT_G$.\vskip 15pt
}
\vskip 10pt
%zzzzzzzzzzzzzzzzzzzzzzzzzzzzzzzzzzzzz Секция 4

{\bf 6.2. Between  $\overline{ K(G)}$  and  $\overline{ G ^*G ^*}$}

% Let $(g_n)_{n\in\omega} $ be an injective sequence  in a group $G$. The set
% $$\{   g_{i_1 }  g_{i_2 } \dots g_{i_n } : 0\leq i_1< i_2< \dots < g_{i_n } < \omega \}  $$ is called an {\it $FP$-set }.

% Given a sequence $(b_n)_{n\in\omega} $  in $G$, we say that the set
% $$\{   g_{i_1 }  g_{i_2 } \dots g_{i_n }  b_{i_n }  : 0\leq i_1< i_2< \dots < i_n  < \omega \}  $$ is a {\it piecewise shifted  $FP$-set }.

\vskip 7pt

%zzzzzzzzzzzzzzzzzzzzzzzzzzzzzzzzzzzzz
{\bf Theorem 6.8.}
{\it  For a group $G$, the following statements hold:
\vskip 10pt

$(i)$ $Sc_G^{\wedge} = cl \{\epsilon p:  \epsilon\in G ^*  , \  p\in \beta G,   \   \epsilon\epsilon=\epsilon\}$;

\vskip 10pt

$(ii)$  $Sc_G^{\wedge}$ is an ideal in $\beta G$ and $p\in Sc_G^{\wedge}$
 if and only if each member of $p$ contains a piecewise shifted  $FP$-set;

\vskip 10pt

$(iii)$   $Sc_G^{\wedge}$ is  the minimal closed ideal in $\beta G$  containing all idempotents of
$G ^*$.\vskip 15pt}

For a group $G$, we put $I _{G,n} = G^{\ast} $,  $I _{G,n+1} = \overline{G^{\ast} I_{G,n}} $  and note that $I _{G,n}$ is an ideal in $\beta G$.
\vskip 15pt

%zzzzzzzzzzzzzzzzzzzzzzzzzzzzzzzzzzzz
{\bf Theorem 6.9.}
{\it  For every group $G$ and $n\in \omega$, we  have \vskip 5pt

$(i) \ $   $I_{G,n+1} \subset I_{G,n }$
\vskip 10pt

$(ii) \ $   $Sc_{G}^{\wedge} \subset I_{G,n }$.}

\vskip 15pt

For a natural number $n$, we denote by  $(G ^*)^n$
the product of $n$ copies of $n$.
Clearly, $\overline{ (G ^*) ^{n+1}}\subseteq   \overline{ (G ^*) ^{n}}$.
and  $ \overline{ (G ^*) ^{n}}\subseteq I_{G,n }$.
\vskip 10pt

{\bf Theorem 6.10.}
{\it  For every  group $G$  and $n\in\omega$, we have   \vskip 5pt

$(i) \ $   $\overline{ (G ^*) ^{n+1}} \subset  \overline{ (G ^*) ^{n}}$;
\vskip 5pt

$(ii) \ $   $Sc_{G}^{\wedge} \subset  \overline{ (G ^*) ^{n}}$.}

\vskip 10pt

{\it Comments.} This section is an extract from \cite{b25}.

%zzzzzzzzzzzzzzzzzzzzzzzzzzzzzzzzzzzz Секция 7

\section{ The combinatorial derivation}

Let $G$  be a group with the identity $e$.
For a subset $A$ of $G$, we denote
$$\triangle (A)= \{ g\in G: |gA\bigcap  A = \infty | \}, $$
observe that
$(\triangle (A))^{-1}= \triangle (A)$, $\triangle (A)\subseteq AA^{-1}$,  and say that the mapping
$$\bigtriangleup : \mathcal{P}_{G}\longrightarrow \mathcal{P}_{G},   \   \  A\longmapsto  \triangle (A)$$
   is the {\it combinatorial derivation.}
   \vskip 10pt

{\bf  Theorem 7.1.}  {\it
For an infinite group $G$ and a subset $A$  of $G,$   the following statements hold

(1)	 $A$ is finite if and only if $\bigtriangleup(A)=\emptyset$;

(2)	  $\bigtriangleup(A)=\{ e \}$ if and only if $A$  is infinite and thin;

(3)	 if  $A$ is thick then $\triangle (A)=G$;

(4)	 if  $A$ is prethick then $\triangle(A)$ is large;}

\vskip 10pt

{\bf  Theorem 7.2.}  {\it Every  infinite group $G$ contains a subset $A$
such that  $G=AA^{-1}$   and $\triangle(A)= \{ e\}$. }
\vskip 10pt

{\bf  Theorem 7.3.}  {\it Let $A$  be a subset of an infinite group $G$ such that $ A=A^{-1}$.   Then there exist two thin  subsets $X$, $Y$  of $G$  such that $\triangle(X\bigcup Y)=A$.}
\vskip 10pt

We consider also the inverse to $\triangle$, multivalued mapping $\nabla$ defined by
 $$\nabla(A)= \{B\subseteq  G: \triangle (B)= A\}.$$

For a family $F$ of subsets of a group $G$, we say that $\mathcal{F}$ is $\triangle$-complete ($\nabla$-complete)  if  $\triangle(A)\in \mathcal{F}$  $(\nabla(A)\subseteq \mathcal{F})$  for each $A\in \mathcal{F}$.
\vskip 10pt

{\bf  Theorem 7.4.}  {\it For every  infinite group $G$, the following statements hold
\vskip 5pt

(1)	 the families of all small and sparse subsets of $G$ is $\nabla$-complete;\vskip 5pt

(2)	 if an ideal $\mathcal{I}$ in $\mathcal{P}_{G}$ is $\triangle$-complete and  $\nabla$-complete  then  $\mathcal{I}=\mathcal{P}_{G}$;\vskip 5pt

(3)	 If  $\mathcal{I}$ is a group ideal in $\mathcal{P}_{G}$, $\mathcal{I} \neq \mathcal{P}_{G}$, then  $\mathcal{I}$  is $\triangle$-complete and  $\mathcal{I}$  is contained in the ideal of all small subsets of} $G$.
\vskip 10pt

{\it Comments.}  More information on combinatorial  derivation in \cite{b26}, \cite{b27}, \cite{b28}.
In particular, Theorem 6.2 from \cite{b26} shows that the trajectory $A\longrightarrow \triangle(A)\longrightarrow \triangle ^{2}(A)\longrightarrow \ldots $ of a subset $A$ of $G$  could be surprisingly complicated: stabilizing, increasing, decreasing, periodic or chaotic. Also \cite{b26} contains some parallels between the combinatorial and topological derivations.

\vspace{3 mm}
CONTACT INFORMATION


\begin{thebibliography}{30}


\bibitem{b1} {I.Protasov,} {\em Selective survey on subset combinatorics of groups}// J. Math. Sciences, {\bf 174} (2011), 486--514.


\bibitem{b2} {I. Protasov, S. Slobodianiuk,} {\em On the subset combinatorics of $G$-spaces}// Algebra Discrete Math.,
  {\bf 17}(2011), 98--109.

\bibitem{b3} {I. Protasov, S. Slobodianiuk,} {\em Partitions of groups}// Math. Stud., {\bf 42}(2014), 115--128.

\bibitem{b4}  {T. Banakh, I. Protasov, S. Slobodianiuk,}  {\em Densities, submeasures and partitions of groups} // Algebra Discrete Math., {\bf 17}(2014), 193--221.


\bibitem{b5}{N. Hindman, D. Strauss,} {\em Algebra in the Stone-\v Cech compactification: theory and applications}, Berlin, New York: Walter de Gruyter, 1998.

\bibitem{b6}  {I. Protasov, S. Slobodianiuk,}  {\em Ultracompanions of subsets of a group} // Comment. Math. Univ. Carolin., {\bf 55} (2014), 257--265.

\bibitem{b7}  {T. Banakh, I. Protasov, S. Slobodianiuk,}  {\em Scattered subsets of groups} // Ukr. Math. J. {\bf 67}:3 (2015) 347--356.


\bibitem{b8} {Ie. Lutsenko, I. Protasov, }  {\em Sparse, thin and other subsets of groups} // Intern. J. Algebra Comp. {\bf 19} (2009) 491--510.

\bibitem{b9}	{A. Kechris,} {\em Classical Descriptive Set Theory}, Springer, 1995.

\bibitem{b10}  {T. Banakh, N. Lyaskovska,} {\em Weakly $P$-small not $P$-small subsets in groups} // Intern. J. Algebra Comput.  {\bf 19} (2008) 1--6.

\bibitem{b11}  {I. Protasov, K. Protasova,} {\em Around $P$-small subsets of groups} // Carpath. Math. Publ. {\bf 6} (2014) 337--341.

\bibitem{b12} {T. Banakh, N. Lyaskovska,} {\em On thin-complete ideals of subsets of groups} // Ukr. Math. J.  {\bf 63}:6 (2011)P 216--225.

\bibitem{b13} {M. Filali, Ie. Lutsenko, I. Protasov,} {\em Boolean group ideals and the ideal structure of $\beta G$} // Math. Stud. {\bf 30} (2008) 1--10.

\bibitem{b14}  {T. Banakh, I. Protasov, K. Protasova,} {\em Descriptive complexity of the sizes of subsets of groups} //Ukr. Mat. J. {\bf 69} (2017), № 9, 1280--1283.

\bibitem{b15} {I. Protasov, S. Slobodianiuk,}  {\em The dynamical look at the subsets of a group} // Appl. Gen. Topol.  {\bf 16}:2 (2015) 217--224.


\bibitem{b16}  { H. Furstenberg,} {\em  Poincare recurrence and number theory} // Bull. Amer. Math. Soc.,
{\bf 5} (1981), № 3, 211–-234.

 \bibitem{b17}  { N. Hindman,} {\em  Ultrafilters and combinatorial number theory} // Lecture Notes in Math., {\bf 571}(1979), 119–-184.

\bibitem{b18}  {V. Bergelson, N. Hindman,} {\em  Quotient sets and density  recurrent sets }// Trans. Amer. Math. Soc., {\bf 364}(2012), 4495–-4531.

\bibitem{b19}  {I.Protasov, }  {\em Filters and topologies on  groups}// Math. Stud., {\bf 3} (1994), 15--28.

\bibitem{b20}  {I. Protasov, K. Protasova,}  {\em On recurrence in $G$-spaces} // Algebra Discrete Math., {\bf 23}:2(2017),
 80--85.

\bibitem{b21}  {T. Banakh, I. Protasov, K. Protasova,}  {\em Ramsey-product subsets of a group} // Math. Stud., {\bf 47} (2017), 145--149.

\bibitem{b22}  { Ie. Lutsenko, I. Protasov,} {\em Thin subsets of balleans}  // Appl. Gen. Topology {\bf 11}
(2010), 89–-93.

\bibitem{b23}{I.Protasov, S. Slobodianiuk,}  {\em Thin subsets of groups} // Ukr. Math. J.   {\bf 65} (2013), 1384--1393.

\bibitem{b24}  {I. Protasov,T. Banakh ,}  {\em Ball Structures and
Colorings of Graphs and Groups}// Math. Stud. Monogr. Ser, {\bf 11}, Lviv: VNTL Publisher, 2003.


\bibitem{b25} {I.Protasov, K. Protasova,} {\it Ideals in $PG$  and $\beta G$}//  preprint, ArXiv: 1704.02494--1.

\bibitem{b26} {I.Protasov,} {\it The combinatorial derivation}//  Appl. Gen. Topology, {\bf 14} (2013), 171--178.

\bibitem{b27} { I.Protasov,} {\it The combinatorial derivation and its inverse mapping}//  Central Europ.J. Math.,
{\bf 11} (2013), 1276--1281.

\bibitem{b28}  {J. Erde, }  {\em A note on combinatorial derivation}// preprint, arxiv: 1210. 7622.



\end{thebibliography}
\end{document}